\documentclass[12pt]{article}

\usepackage{amssymb,amsmath,amsfonts,enumerate}

\title{
\author{\bf{Lech Pasicki}}
\bf{Contractions, inwardness, tool theorems}
\footnote{MSC: 54H25, 54E35, 47H09, 47H10}}

\newtheorem{theorem}{\indent Theorem}[section]
\newtheorem{lemma}[theorem]{\indent Lemma}

\newtheorem{definition}[theorem]{\indent Definition}
\newtheorem{corollary}[theorem]{\indent Corollary}
\newtheorem{remark}[theorem]{\indent Remark}

\newtheorem{problem}[theorem]{\indent Problem}

\newcommand{\mN}{\mbox{$\mathbb{N}$}}
\newcommand{\Api}{\mbox{$\mathcal{A}$}}
\newcommand{\tn}{\mbox{$t_{n}$}}
\newcommand{\xk}{\mbox{$x_{k}$}}
\newcommand{\xo}{\mbox{$x_{0}$}}
\newcommand{\xp}{\mbox{$x_{1}$}}
\newcommand{\xn}{\mbox{$x_{n}$}}
\newcommand{\xnp}{\mbox{$x_{n+1}$}}
\newcommand{\aeps}{\mbox{$\alpha + \epsilon$}}
\newcommand{\xknn}{\mbox{$(x_{k_{n}})_{n \in \mathbb{N}}$}}
\newcommand{\xkn}{\mbox{$x_{k_{n}}$}}
\newcommand{\xnn}{\mbox{$(x_{n})_{n \in \mathbb{N}}$}}
\newcommand{\yn}{\mbox{$y_{n}$}}
\newcommand{\zn}{\mbox{$z_{n}$}}
\newcommand{\An}{\mbox{$A_{n}$}}

\newcommand{\fXX}{\mbox{$f\colon X\to X$}}
\newcommand{\dFc}{\mbox{$d(F(\cdot),\cdot)$}}
\newcommand{\FXX}{\mbox{$F\colon X\to 2^{X}$}}
\newcommand{\FXY}{\mbox{$F\colon X\to 2^{Y}$}}
\newcommand{\Xd}{\mbox{$(X,d)$}}
\newcommand{\Yd}{\mbox{$(Y,d)$}}

\newcommand{\li}{\mbox{$\lim_{n \to \infty}$}}
\newcommand{\lin}{\mbox{$\liminf_{n \to \infty}$}}
\newcommand{\lis}{\mbox{$\limsup_{n \to \infty}$}}
\newcommand{\ej}{\mbox{$\epsilon_{1}$}}
\newcommand{\jep}{\mbox{$1-\epsilon$}}
\newcommand{\dcF}{\mbox{$d(\cdot,F(\cdot))$}}
\newcommand{\dxF}{\mbox{$d(x,F(x))$}}
\newcommand{\dzF}{\mbox{$d(z,F(z))$}}
\newcommand{\dzFx}{\mbox{$d(z,F(x))$}}
\newcommand{\delxz}{\mbox{$\delta(x,z)$}}
\newcommand{\fix}{\mbox{$\varphi(x)$}}
\newcommand{\fixn}{\mbox{$\varphi(x_{n})$}}
\newcommand{\fiz}{\mbox{$\varphi(z)$}}
\newcommand{\pn}{\mbox{$p_{n}$}}
\newcommand{\tIXx}{\mbox{$\tilde{I}_{X}(x)$}}
\newcommand{\xm}{\mbox{$x_{m}$}}
\newcommand{\xXFx}{\mbox{$x \in X \setminus F(x)$}}
\newcommand{\zXx}{\mbox{$z \in X \setminus \{x\}$}}
\newcommand{\cf}{{\it cf}}
\newcommand{\eg}{{\it e.g}}

\begin{document}
\maketitle
\vspace{1 in}

\begin{abstract}
\par The paper is devoted to the fixed point theory in four aspects: of contractions, nonexpansive mappings,
generalized inward mappings, and of the tool theorems. The manuscript was written about ten years ago. 
\par At first Nadler's concept of contraction for multivalued mappings is replaced here by a more
general, and yet elegant condition: {\it for some $\alpha + \epsilon <1$, and each $x \in X$ there exists
a $y \in F(x)$ such that $d(F(y),y) \leq \alpha d(y,x) \leq (\alpha + \epsilon) d(F(x),x)$}.
\par For ``nonexpansive'' mappings we apply bead spaces that are more general than uniformly convex spaces,
and our requirements on mappings are weaker than nonexpansivity in the sense of the Hausdorff distance.
\par In the last, third section the Caristi theorem is replaced by more specialized ``tools'', and we apply 
them to obtain stronger fixed point theorems on generalized inward mappings. In particular, if for each  
$x \in X$ a nearest point of $F(x)$ belongs to the generalized inward set, then the 
values of $F$ need not to be closed.
\end{abstract}

\par The following three theorems are well known:
\begin{theorem}[{\rm Banach\bf}]
\label{Th01}
Let $ \Xd $ be a complete metric space, $\alpha < 1$, and let a mapping $\fXX$ satisfy
\begin{equation}
\label{co1}
 d(f(y),f(x)) \leq \alpha d(y,x), \quad x,y \in X.
\end{equation}
Then $f$ has a fixed point.
\end{theorem}
\par Let $2^Y$ be the family of all subsets of $Y$. We say that $ \FXY $ is a (multivalued)
mapping if $F(x) \neq \emptyset $, $x \in X \neq \emptyset$. The graph of $F$ is the set 
$graph(F) = \{(x,y)\colon x \in X \mbox{, } y \in F(x) \}$.
\par If $ \Yd $ is a metric space, $A,B \subset Y$ are nonempty, and
\[
\max \mbox{\{}\sup\limits_{x \in A} d(x,B), \sup\limits_{y \in B} d(A,y) \mbox{\}}
\]
is finite, then the latter is denoted by $D(A,B)$ (the Hausdorff distance).
\begin{theorem}[{\rm Nadler \cite[Theorem 5]{Na}\bf}]
\label{Th02}
Let $ \Xd $ be a complete metric space, $\alpha < 1$, and let a mapping $\FXX$
with bounded, closed values satisfy
\begin{equation}
\label{co2}
 D(F(y),F(x)) \leq \alpha d(y,x), \quad x,y \in X.
\end{equation}
Then $F$ has a fixed point.
\end{theorem}
\par A mapping satisfying \eqref{co1} or \eqref{co2} for $\alpha < 1$ is called
$\mathbf{\alpha}$-{\bf contraction}, and for $\alpha = 1$ such a mapping is {\bf nonexpansive}.
\begin{theorem}[{\rm Lim \cite[Theorem 1]{Lf}\bf}]
\label{Th03}
Let $X$ be a bounded, closed, convex set in a uniformly convex Banach space, and $\FXX$ a nonexpansive
mapping with compact values. Then $F$ has a fixed point.
\end{theorem}

\section{Contractions}
\label{Con}
\par Let us investigate conditions \eqref{co1}, \eqref{co2}. From \eqref{co1} it follows that
\[
	d(f(y),y) \leq \alpha d(y,x), \quad x\in X \mbox{, } y=f(x).
\]
Similarly, from \eqref{co2} we obtain
\[
	d(F(y),y) \leq D(F(y),F(x)) \leq \alpha d(y,x), \quad x \in X \mbox{, } y \in F(x).
\]
Consequently, for $y \in F(x)$ such that $d(y,x) \leq (1+\epsilon/\alpha)d(F(x),x)$
( $\alpha > 0$) we get
\[
	d(F(y),y) \leq \alpha d(y,x) \leq (\aeps) d(F(x),x), \quad x \in X.
\]
Now, we can see that for any $\epsilon > 0$ it follows from condition \eqref{co2} that
\begin{equation}
\label{co3}
\begin{split}
 &\mbox{for each } x\in X \mbox{ there exists a } y\in F(x) \mbox{ such} \\
 &\mbox{that } d(F(y),y) \leq \alpha d(y,x) \leq (\aeps) d(F(x),x), 
\end{split}
\end{equation}
and the dependence is valid also for $\alpha=0$.
\par Let us recall that a mapping $\FXX$ is an $\mathbf{\alpha}$-{\bf step} \cite[Definition 17]{Pt} if
\begin{equation}
\label{co4}
\begin{split}
   &\mbox{for each } x \in X, y \in F(x) \mbox{ there exists} \\
   &\mbox{a } z \in F(y) \mbox{ such that }  d(z,y) \leq \alpha d(y,x)
\end{split}
\end{equation}
holds.
\par Clearly, from $d(F(y),y) \leq d(z,y)$, $z \in F(y)$ it follows that \eqref{co3} is more general
than \eqref{co4}.
\par These considerations suggest that condition \eqref{co3} is a natural component for
an extension of the Nadler theorem. 
\begin{lemma}
\label{Le11}
Let $ \Xd $ be a metric space, and $\FXX$ a mapping satisfying \eqref{co3} 
with $\aeps < 1$. Then there exists a Cauchy sequence $\xnn$ in $X$ such that
$\xnp \in F(\xn),$ $n \in \mN$.
\end{lemma}
Proof. 
It is sufficient to consider $\alpha > 0$. Let $\xp \in X$
be arbitrary. Then for $x:=\xn$, and $y := \xnp$ as in \eqref{co3}, $n \in \mN$,
we obtain
\begin{equation*}
\begin{split}
 &\alpha d(\xnp,\xn) \leq (\aeps) d(F(\xn),\xn) \leq 
   (\aeps) \alpha d(\xn,x_{n-1}) \leq \\
 &(\aeps)^2 d(F(x_{n-1})),x_{n-1})\leq \ldots
\leq (\aeps)^n d(F(\xp),\xp).
\end{split}
\end{equation*}
Consequently, $d(\xnp,\xn) \leq (\aeps)^n C$, $n \in \mN$, holds, and 
$\sum^{\infty}_{n=1} d(\xnp,\xn)$ converges, which means (as in Banach's proof) that
$\xnn$ is a Cauchy sequence.  $\quad \square$
\par For $\epsilon = 0$ Lemma \ref{Le11} can be reformulated as follows
\begin{lemma}
\label{Le12}
Let $ \Xd $ be a metric space, and $\FXX$ a mapping satisfying 
\begin{equation}
\label{co5}
\begin{split}
 \mbox{for each } x\in X \mbox{ there exists a } y\in F(x) \mbox{ such that} \\
 d(y,x)=d(F(x),x) \mbox{, and } d(F(y),y) \leq \alpha d(F(x),x)
\end{split} 
\end{equation}
for an $\alpha < 1$. Then there exists a Cauchy sequence $\xnn$ in $X$ such that
$\xnp \in F(\xn),$ $n \in \mN$.
\end{lemma}
\par Clearly, if $F(x)$ is compact (or closed convex in a Hilbert space), then it contains
a point nearest to $x$.
\par In connection with Lemmas \ref{Le11}, \ref{Le12} we are interested in the following dependence:
\begin{equation}
\label{co6}
\begin{split}
 &\mbox{for each } \xnn \mbox{ in } X \mbox{ that converges, say to an } x \in X \mbox{, } \\
 &\mbox{if }  \li d(F(\xn),\xn)=0 \mbox{, then } d(F(x),x)=0. 
\end{split}
\end{equation}
It holds, \eg., if $\varphi\colon X\to [0,\infty)$ defined by $\varphi(x)=d(F(x),x)$, $x \in X$, is 
continuous in its zeros. This happens in particular for lower semicontinuous $\varphi$.
This case requires closer attention.
\par Let us recall that for topological spaces $X,Y$ a mapping $ \FXY $ is {\bf usc}
if for every $x \in X$, and any neighborhood $V$ of $F(x)$ there exists a
neighborhood $U$ of $x$ such that $F(U) \subset V$.
\par The subsequent lemma extends \cite[Lemma 2]{Do}:
\begin{lemma}
\label{Le13}
Let $X$ be a set in a metric space $\Yd$, and $\FXY$ a mapping that is usc or continuous
with respect to the Hausdorff metric in $Y$. Then $\dFc$ is lower semicontinuous.
\end{lemma}
Proof. 
In view of the continuity assumption on $F$, if $\xnn$ converges in $X$ to a point $x$, then
for $\yn \in F(\xn)$ there exist $\zn \in F(x)$
satisfying $\li d(\zn,\yn) = 0 $. Hence, for $(y_{n})_{n \in \mathbb{N}}$ such that 
$\lin d(\yn,\xn) = \lin d(F(\xn),\xn)$ we obtain
\begin{equation*}
\begin{split}
&d(F(x),x) \leq \lin d(\zn,x) \leq \lin [d(\zn,\yn)+ \\
&d(\yn,\xn)+d(\xn,x)] = \lin d(F(\xn),\xn).\qquad \square
\end{split}
\end{equation*}
\par Another way of obtaining \eqref{co6} is to assume that the graph of  $F$ is closed. Clearly, all values of
such a mapping must be closed. If $Y$ is a regular space, $ \FXY $ is usc, and all values of $F$ are closed,
then the graph of $F$ is closed (see the proof of \cite[Lemma, p. 285]{Br}). The same holds for 
$(Y,d)$ being a metric space, and $ \FXY $ a closed-valued mapping that is continuous
with respect to the Hausdorff metric in $2^Y$. Therefore, it is better to assume that the graph of $F$ is 
closed.
\par In view of the considerations preceding Lemma \ref{Le11}, the next theorem extends the
theorem of Nadler, and \cite[Theorem 21]{Pt}. It seems to be well placed between the Banach theorem
and the Nadler one.
\begin{theorem}
\label{Th14}
Let $ \Xd $ be a complete metric space, and $\FXX$ a closed-valued mapping satisfying \eqref{co6}
(e.g. if the graph of $F$ is closed).
If \eqref{co3} for some $\aeps < 1$ is satisfied, then $F$ has a fixed point.
\end{theorem}
Proof. 
Let $\xnn$ be a Cauchy sequence as in Lemma \ref{Le11}. Consequently, $\xnn$ converges, say to $x$,
$\Xd$ being complete. Then $d(F(x),x)=0$ (see \eqref{co6}), and $x \in F(x)$ this last being closed.
$\quad \square$
\par Theorem \ref{Th14} is more general than \cite[Theorem 2]{Do}, as we do not assume $\dFc$ to be lower
semicontinuous.
\par The subsequent theorem is more elegant but weaker:
\begin{theorem}
\label{Th15}
Let $ \Xd $ be a complete metric space, and $\FXX$ a mapping satisfying \eqref{co6}. If \eqref{co5} for 
an $\alpha < 1\: $ holds, then $F$ has a fixed point.
\end{theorem}
Proof. 
For $\xnn$ as in Lemma \ref{Le12} there exists an $x \in X$ such that $d(F(x),x)=0$, 
and in view of \eqref{co5} $x \in F(x)$. $\quad \square$
\par For a relaxed version of condition \eqref{co5} we obtain
\begin{lemma}
\label{Le16}
Let $ \Xd $ be a metric space, and $\FXX$ a mapping satisfying 
\begin{equation}
\label{co7}
\begin{split}
 &\mbox{for each } x\in X \mbox{ there exists a } y\in F(x) \\
 &\mbox{such that } d(F(y),y) \leq \alpha d(F(x),x) 
\end{split}
\end{equation}
for an $\alpha < 1 $. Then there exist $\xn \in X$ such that $\li d(F(\xn),\xn) = 0$.
\end{lemma}
Proof. 
For arbitrary $\xp \in X$, and $x:=\xn $, $y := \xnp$ as in \eqref{co7}, 
$n \in \mN$, we obtain
\[
	d(F(\xnp),\xnp) \leq \alpha d(F(\xn),\xn) \leq \cdots \leq \alpha^n d(F(\xp),\xp).\qquad \square
\]
\par Clearly, condition \eqref{co7} does not guarantee $\xnn$ to be a Cauchy sequence. We have
the following:
\begin{theorem}
\label{Th17}
Let $ \Xd $ be a complete metric space, and $\FXX$ a closed-valued mapping satisfying \eqref{co6}
(e.g. if the graph of $F$ is closed). If \eqref{co7} holds for some $\alpha < 1$, and $\overline{F(X)}$
is compact, then $F$ has a fixed point.
\end{theorem}
Proof. 
In view of Lemma \ref{Le16} there exist a sequence $\xnn$, and $\yn \in F(\xn)$ such that
$\li d(\xn,\yn) = 0$ . The relative compactness of $\{\yn\colon n \in \mN \}$ yields the relative
compactness of $\{\xn\colon n \in \mN \}$. Therefore, $\xnn$ has a convergent subsequence, say $\xknn$.
Then for  $x = \li \xkn = \li y_{k_{n}}$ we have $x \in F(x)$  (see \eqref{co6}). $\quad \square$

\section{Nonexpansive mappings}
\label{Non}
\par Lim's theorem (Theorem \ref{Th03}) concerns uniformly convex spaces. Let us present a similar idea for metric
spaces.
\begin{definition}[{\rm \cite[Definition 6]{Ps}\bf}]
\label{De21}
  A metric space $\Xd$ is a {\bf bead space} if the following is valid: 
\begin{equation}
\label{co8}
\begin{split}
&\mbox{for every } r > 0, \beta > 0 \mbox{ there exists a } \delta > 0 \mbox{ such that} \\
&\mbox{for each } x,y \in X \mbox{ with } d(x,y) \geq \beta \mbox{ there exists a } z \in X \\
&\mbox{such that } B(x,r+\delta) \cap B(y,r+\delta) \subset B(z,r - \delta).
\end{split}
\end{equation}
\end{definition}
\par In particular, each convex set in a uniformly convex space is a bead space
\cite[Example 3]{Pf} with $z=(x+y)/2$  in condition \eqref{co8} (\cf. \cite[Theorem 4]{Pu}).
\par A normed space is a bead space if and only if it is uniformly convex \cite[Theorem 14]{Pu}.
On the other hand, there exist bead spaces which are not convex sets in uniformly convex spaces
\cite[Example 3]{Pt}.
\begin{definition}[{\rm \cite[Definition 11]{Ps}\bf}]
\label{De22}
  Let $ \Xd $ be a metric space, and $\Api$ a  family of nonempty bounded
subsets of $X$.  An $x \in X$ is a {\bf central point} for $\Api$ if 
\begin{equation}
\label{co9}
\begin{split}
  &r(\Api) := inf\{t \in (0,\infty): \mbox{ there exist } A \in \Api, z \in X
      \mbox{ such that } A \subset \\
  &B(z,t) \} =inf \{t \in (0,\infty):
      \mbox{ there exists } A \in \Api \mbox{ with } A \subset  B(x,t)\}.
\end{split}
\end{equation}
The {\bf centre} $c(\Api)$ for $\Api$ is the set of all central points for $\Api$, and
$r(\Api)$ is the {\bf radius} of $\Api$.
\end{definition}
It should be noted that $r(\Api)$ is defined by condition \eqref{co9} even if
$c(\Api)=\emptyset$.
\par If $\Api$ is a family of nonempty, and bounded sets directed by $\supset$ in a bead space, then
$c(\Api)$ consists of at most one point, and it is a singleton if the bead space under
consideration is complete (\cf. \cite[Lemma 12]{Ps}).
\par If $\xnn$ is a bounded sequence of points of $X$, then for $\An =
\{\xk: k \geq n \}$, and $\Api = \{\An: n \in \mN\}$, $c(\xnn):=c(\Api)$,
$r(\xnn):=r(\Api)$ are, respectively, the {\bf (asymptotic) centre}, and the {\bf (asymptotic)
radius} of $\xnn$ (see \cite{Ed}).
\par From \cite[Lemma 12]{Ps} it follows that the centre of any bounded sequence in a bead space
is at most a singleton, and that it is nonempty if such a space is complete (\cf. \cite{Pt}).
\par Let us recall (see \cite{Go}) that a bounded sequence $\xnn$ in a metric space $\Xd$ is 
{\bf regular} if $r(\xknn)=r(\xnn)$ holds for each of its subsequences $\xknn$; a regular sequence
$\xnn$ is {\bf almost convergent} if $c(\xknn)=c(\xnn) $. It is known that any bounded sequence in
a metric space contains a regular subsequence \cite[Lemma 2]{Go}. On the other hand, any regular sequence in
a bead space is almost convergent \cite[Lemma 13]{Pt}.
\par Now, we are ready to prove the following:
\begin{theorem}
\label{Th23}
Let $\Xd$ be a bead space, and $\FXX$ a mapping. Assume that $\xnn$ is a regular
sequence, $x \in c(\xnn)$, $\yn \in F(\xn)$ are such that $\li d(\yn,\xn) = 0$, and
\begin{equation}
\label{co10}
 \lis d(F(x),\yn) \leq  \lis d(x,\xn)
\end{equation}
holds. If $F(x)$ is compact, then $x \in F(x)$.
\end{theorem}
Proof. 
Let us adopt $r=r((\xnn))$. There exist $\zn \in F(x) $ such that
\begin{equation*}
\begin{split}
&\lis d(\zn,\xn) \leq \lis [d(\zn,\yn) + d(\yn,\xn)] =\\
& \lis d(F(x),\yn) \leq \lis d(x,\xn) = r 
\end{split}
\end{equation*}
(see \eqref{co10}). The sequence
$(z_{n})_{n \in \mathbb{N}}$ has a subsequence $(z_{k_{n}})_{n \in \mathbb{N}}$ that converges, say to a point 
$z$, in the compact set $F(x)$. Hence, $\lis d(z,\xkn) \leq r$, and in view of \cite[Lemma 13]{Pt}
$z \in c(\xknn) = c(\xnn) = \{x\}$ as $\xnn$ is regular. Thus, we have
$x = z \in F(x)$.  $\quad \square$
\par Condition \eqref{co10} is satisfied, \eg., if
\begin{equation}
\label{co11}
 \yn \in F(\xn) \mbox{ are such that }
 d(F(x),\yn) \leq  d(x,\xn), \quad n \in \mN
\end{equation}
holds. In turn, \eqref{co11} follows from
\begin{equation}
\label{co12}
   F(\xn) \subset \overline{B}(F(x),d(x,\xn)), \quad n \in \mN,
\end{equation}
and \eqref{co12} is clearly satisfied for nonexpansive $F$.
\par Let us recall once again that any bounded sequence in a metric space has a regular subsequence, and
the centre of any bounded sequence in a complete bead space is a singleton.
\par By putting these facts together we obtain the following simplified version of Theorem \ref{Th23}:
\begin{theorem}
\label{Th24}
Let $\Xd$ be a complete bead space, and $\FXX$ a bounded mapping satisfying \eqref{co7} for an
$\alpha <1$. Assume that a regular sequence $\xnn$, and $\yn \in F(\xn)$ are such that
$\li d(\yn,\xn) = 0$ (in view of Lemma \ref{Le16} such $\yn$ exist), and for $x \in c(\xnn)$
\eqref{co11} holds. If $F(x)$ is compact, then $x \in F(x)$.
\end{theorem}
\par The formulation is clearer in the case of condition \eqref{co12}:
\begin{theorem}
\label{Th25}
Let $\Xd$ be a complete bead space, and $\FXX$ a bounded mapping satisfying \eqref{co7} for an
$\alpha <1$. If for the centre $\{x\}$ of any regular sequence $\xnn$ such that
$\li d(F(\xn),\xn) = 0$  the set $F(x)$ is compact, and \eqref{co12} holds, then
$F$ has a fixed point.
\end{theorem}
The theorem also follows from \cite[Theorem 14]{Pt}.
\par If we assume that all values of $F$ are compact, then Theorem \ref{Th25} yields 
\begin{theorem}
\label{Th26}
Let $\Xd$ be a complete bead space, and $\FXX$ a bounded, compact-valued mapping satisfying \eqref{co7}
for an $\alpha <1$. If for the centre $\{x\}$ of any regular sequence $\xnn$ such that
$\li d(F(\xn),\xn) = 0$  \eqref{co12} holds, then $F$ has a fixed point.
\end{theorem}
\par The next theorem is a consequence of Theorem \ref{Th25}, and it has a ``classical'' outlook.
\begin{theorem}
\label{Th27}
Let $\Xd$ be a complete bead space, and $\FXX$ a bounded, compact-valued nonexpansive mapping. If
there exists a sequence $\xnn$ such that $\li d(F(\xn),\xn) = 0$ (e.g. \eqref{co7} holds for
an $\alpha <1 $), then $F$ has a fixed point.
\end{theorem}
\par Theorem \ref{Th27} is close to Lim's Theorem \ref{Th03}, and for condition \eqref{co7} it is situated
between the Nadler Theorem \ref{Th02}, and the Lim theorem.
\begin{remark}
\label{Re28}
If for each $x \in X$ there exists a $y \in F(x)$ such that $d(y,x) = d(F(x),x)$, 
then we define a selection $g\colon X\to X$ for $F$ by taking $g(x)=y$. Our paper \cite{Pt} contains some 
theorems concerning the case of (multivalued) selections. 
\end{remark}
\begin{problem}
\label{Pr29}
In his original proof Lim has shown that for any bounded nonexpansive mapping in a uniformly convex
space $X$ there exist  $\xn \in X$ with $\li d(F(\xn),\xn) = 0$. If this was true for the bead spaces, 
condition \eqref{co7} could be disregarded. A proof or a counterexample is needed.
\end{problem}
\section{Tool theorems, and inward mappings}
\label{Too}
\par Let us recall the Caristi theorem.
\begin{theorem}[{\rm Caristi \cite[Theorem (2.1)']{Ca}\bf}]
\label{Th31}
Let $\Xd$ be a complete metric space, and $g \colon X \to X$ a mapping. If $\varphi \colon X \to [0,\infty)$ is a 
lower semicontinuous mapping such that
\[
d(x,g(x)) \leq \varphi(x) - \varphi(g(x)), \quad x \in X,
\]
then $g$ has a fixed point.
\end{theorem}
\par Our first ``tool'' theorem is the following:
\begin{theorem}
\label{Th32}
Let $(X,\delta)$ be a complete metric space, and $X \subset Y$. Assume that $\varphi\colon X\to \mathbb{R}$ is 
a lower semicontinuous mapping having a finite lower bound. If $ \FXY $ is a mapping satisfying
\begin{equation}
\label{co13}
\begin{split}
 &\mbox{for each  } \xXFx \mbox{  there exists a  } \zXx \\
 &\mbox{such that } \delxz \leq  \fix - \fiz, 
\end{split}
\end{equation}
then $F$ has a fixed point.
\end{theorem}
Proof. 
We apply the idea from \cite{Pc}. In view of the Teichm\"uller-Tukey lemma there exists an 
$\xo$ belonging to a maximal set $A \subset X$ such that all $x,z \in A$ satisfy
\[
  \delxz \leq \mid \fix - \fiz \mid. 
\]
Let us adopt $\gamma = \inf \{ \fiz\colon z \in A\}$, and suppose that $\gamma < \fiz $,  
$z \in A $. Then there exists a sequence $\xnn$ such that $(\varphi(x_{n}))_{n \in \mathbb{N}}$ decreases 
to $\gamma$. For each $m < n$ we have
\begin{equation*}
\begin{split}
 &\delta(\xm,\xn) \leq \delta(\xm,x_{m+1})+ \cdots + \delta(x_{n-1},\xn) \leq \\
 &\varphi(x_{m}) - \varphi(x_{m+1}) + \cdots +\varphi(x_{n-1}) - \fixn = \varphi(x_{m}) - \fixn,
\end{split}
\end{equation*}
and therefore $\xnn$ is a Cauchy sequence converging, say to an $x \in X$ ($X$ is complete).
For any $z \in A$, and large $n$ we have $\delta(x_{n},z) \leq \fiz - \fixn$,
\[
 \delxz \leq \delta(x,x_{n}) + \delta(x_{n},z) \leq \delta(x,x_{n})+ \fiz - \fixn \leq
 \delta(x,x_{n}) + \fiz - \fix,
\]
($\varphi$ is lower semicontinuous), and consequently, $\delxz \leq \fiz - \fix$, {\it i.e.} $x \in A$. 
In addition, $\gamma \leq \fix < \fixn$ means $\fix = \gamma$. Suppose $x \notin F(x)$, and 
let $z \neq x$ be as in \eqref{co13}, {\it i.e.} 
\[
 \fiz \leq \fix - \delxz < \fix = \gamma.
\]
On the other hand, for any $y \in A$ we have
\[
 \delta(y,z) \leq \delta(y,x) + \delxz \leq \varphi(y) - \fix + \fix - \fiz = \varphi(y) - \fiz,
\]
{\it i.e.} $z \in A$, and therefore, $\gamma \leq \fiz$ which contradicts  $\fiz < \gamma $.  $\quad \square$
\par The direct proof of the preceding general theorem is short, and the theorem itself is a convenient tool. 
On the other hand, Theorems \ref{Th31}, \ref{Th32} are equivalent as shown by the 
following reasoning: \par Theorem \ref{Th31} is a consequence of Theorem \ref{Th32} for $F(x) := \{g(x)\}$, 
and $z:=g(x)$. On the contrary, assume $F$ as in Theorem \ref{Th32} has no fixed point. Then for 
$g(x) = z$ ($z$ as in \eqref{co13}), $x \in X$, $g\colon X\to X$ is a fixed point free mapping satisfying 
the assumptions of Theorem \ref{Th31}.
\par The method presented in the proof of Theorem \ref{Th32} is used, in a modified way, to prove 
theorems \ref{Th38}, \ref{Th312}.
\par It is worth noting that \cite[Theorem 23]{p30} is a far extension of Theorem \ref{Th32}, and maybe it could 
a be starting point to modernize the present paper.
\par As a consequence of Theorem \ref{Th32} we obtain the following:
\begin{theorem}
\label{Th33}
Let $\Yd$ be a metric space, $X \subset Y$, and let $(X,\delta)$ be a complete metric space. If $\FXY$ 
is a mapping such that $\dcF$ is lower semicontinuous, and 
\begin{equation*}
\begin{split}
 &\mbox{for each  } \xXFx \mbox{  there exists a  } \zXx \\
 &\mbox{such that } \delxz \leq  \dxF - \dzF 
\end{split}
\end{equation*}
holds, then $F$ has a fixed point.
\end{theorem}
Proof. 
We apply Theorem \ref{Th32} to $\fix = \dxF$, $x \in X$.  $\quad \square$
\par Theorem \ref{Th33} is also a far consequence of \cite[Theorem 24]{Pv} which, alas, has a much longer proof.
\par Theorem \ref{Th33} applies, \eg., in proving the following:
\begin{theorem}
\label{Th34}
Let $X$ be a complete set in a metric space $(Y,d)$. If $ \FXY $ is an $\alpha$-contraction,
$\aeps < 1$, and
\begin{equation}
\label{co14}
\begin{split}
 &\mbox{for each  } \xXFx \mbox{  there exists a  } \zXx \\
 &\mbox{such that } (\jep) d(x,z) \leq \dxF - \dzFx
\end{split}
\end{equation}
is satisfied, then 
\begin{equation}
\label{co15}
\begin{split}
 &\mbox{for each  } \xXFx \mbox{  there exists a  } \zXx \\
 &\mbox{such that } (1-\alpha-\epsilon) d(x,z) \leq \dxF - \dzF
\end{split}
\end{equation}
holds, and consequently, $F$ has a fixed point.
\end{theorem}
Proof. 
We have
\begin{equation*}
\begin{split}
 &\dzF \leq d(z,F(x)) + D(F(x),F(z)) \leq \dzFx + \alpha d(x,z) = \\
 &\dzFx - \dxF + \alpha d(x,z) + \dxF, \quad x,z \in X.
\end{split}
\end{equation*}
Consequently, in view of \eqref{co14} for $\xXFx$ there exists a $\zXx$ such that
\[
 \dzF \leq \alpha d(x,z) + \dxF - (\jep) d(x,z),
\]
{\it i.e.} \eqref{co15} holds, and we can apply Theorem \ref{Th33} for $\delta = (1-\alpha-\epsilon)d$ 
(see Lemma \ref{Le13}). $\quad \square$
\par Condition \eqref{co14} has much in common with generalized inwardness. Let $\Yd$ be a metric space.
Then for $x,t \in Y$ a metric segment is defined by
\[
 [x,t]=\{s \in Y\colon d(x,s)+d(s,t)=d(x,t)\},
\]
and $(x,t]=[x,t] \setminus \{x\}$. For a nonempty set $X \subset Y$ a generalized inward set
$\tIXx$ can be defined by (\cf. \cite[p. 1209]{Ma})
\begin{equation*}
\begin{split}
 &\tIXx = \{t \in Y\colon \mbox{ for each } \beta > 0 \mbox{ there exists an} \\
 &s \in (x,t] \mbox{ such that } d(z,s) \leq \beta d(x,s)\}.
\end{split}
\end{equation*}
The inward set $I_{X}(x)$ in a normed space is defined by
\[
 I_{X}(x) = x + \{\lambda (z-x)\colon z \in X \mbox{, } \lambda \geq 1\},
\]
and the weakly inward set is $\overline{I_{X}(x)}$.

Maciejewski has proved \cite[Lemma 1.3]{Ma} that if $Y$ is a normed space, then $\overline{I_{X}(x)} \subset \tIXx$ .
\begin{lemma}
\label{Le35}
Let $\Yd$ be a metric space, $X,C$ nonempty subsets of $Y$, and $x \in X \subset C$. Then for 
each $\epsilon > 0 $, $t \in C \cap \tIXx$ there exists a $\zXx$ such that
\[
 (\jep) d(x,z) \leq d(x,t) - d(z,t).
\]
If $t$ is the nearest point of $C$ to $x$, and $t \in \tIXx $, then
\[
 (\jep) d(x,z) \leq d(x,C) - d(z,C).
\]
\end{lemma}
Proof. 
For a $t \in C$ let $s \in (x,t]$, $z \in X$ be such that $d(z,s) \leq \beta d(x,s)$.
We have
\begin{equation*}
\begin{split}
 &d(x,z) + d(z,t) - d(x,t) \leq d(x,s) + d(s,z) + d(z,s) + d(s,t) - d(x,t) \leq \\
 &2 d(z,s) + d(x,s) + d(s,t) - d(x,t) = 2 d(z,s) \leq 2 \beta d(x,s).
\end{split}
\end{equation*}
From
\[
 d(x,s) \leq d(x,z) + d(z,s) \leq d(x,z) + \beta d(x,s)
\]
we obtain $d(x,s) \leq d(x,z)/(1-\beta)$ for $\beta < 1$. Hence it follows that
\[
 d(x,z) + d(z,t) - d(x,t) \leq 2 \beta d(x,z)/(1-\beta) \leq \epsilon d(x,z),
\]
{\it i.e.}
\[
 (\jep) d(x,z) \leq d(x,t) - d(z,t)
\]
for small $\beta > 0$. If $t$ is the nearest point to $x$ in $C$, then we have
\[
 (\jep) d(x,z) \leq d(x,C) - d(z,t) \leq d(x,C) - d(z,C). \qquad \square
\]
\par From Theorem \ref{Th34}, and Lemma \ref{Le35} when applied to $C = F(x)$ we obtain the following
extension of \cite[Theorem 2.3]{Ma} (we do not demand $F(x)$ to be closed), and of \cite[Theorem 3.3]{Xu}:
\begin{theorem}
\label{Th36}
Let $X$ be a complete set in a metric space $(Y,d)$. If $\FXY$ is an $\alpha$-contraction, each 
$x \in X$ has a nearest point in $F(x)$, and such a point belongs to $\tIXx$, then
for each $\epsilon > 0$ \eqref{co15} is satisfied, and consequently, $F$ has a fixed point.
\end{theorem}
\par For the inwardness condition we have a ``strict'' result:
\begin{theorem}
\label{Th37}
Let $X$ be a complete set in a normed space $Y$. If $F\colon X \to 2^Y$ is an $\alpha$-contraction, each $x \in X$ 
has a nearest point in $F(x)$, and such a point belongs to $I_{X}(x)$, then
\eqref{co15} is satisfied for $\epsilon = 0$, and $F$ has a fixed point.
\end{theorem}
Proof. 
For $y = x + \lambda(z-x)$ we have $d(y,z) = d(y,x) - d(z,x)$. Hence it follows (see \eqref{co2}) that
\begin{equation*}
\begin{split}
 &d(F(z),z) \leq D(F(z),F(x)) + d(F(x),z) \leq \\
 &\alpha d(z,x) + d(y,z) = \alpha d(z,x) + d(y,x) - d(z,x),
\end{split}
\end{equation*}
and for $d(y,x) = d(F(x),x)$ we get \eqref{co15} with $\epsilon = 0$. Now, we can apply 
Theorem \ref{Th33} for $\delta = (1-\alpha)d$ (see Lemma \ref{Le13}). $\quad \square$
\par Let us present another ``tool'' theorem.
\begin{theorem}
\label{Th38}
Let $X$ be a set in a complete metric space $(Y,d)$. Assume that a mapping $\FXY$ has 
closed graph, and $\delta$ is a metric on $X \cup F(X)$ equivalent to $d$. If 
\begin{equation}
\label{co16}
\begin{split}
 &\mbox{for each  } \xXFx \mbox{, } t \in F(x) \mbox{ there exist } \zXx \mbox{, } v \in F(x) \\
 &\mbox{such that } \delxz \leq d(x,t) - d(z,v) \mbox{, and  } \delta(t,v) \leq k(d(x,t) - d(z,v))
\end{split}
\end{equation}
holds for a $k > 0 $, then $F$ has a fixed point.
\end{theorem}
Proof. 
Assume $\xo \notin F(\xo)$, and let $A \neq \emptyset$ be a maximal set in $G = graph(F)$ 
such that for $(x,t), (z,v) \in A$
\[
 \delxz \leq \mid d(x,t) - d(z,v) \mid, \quad \delta(t,v) \leq k \mid d(x,t) - d(z,v) \mid 
\]
are satisfied. Let us adopt $\gamma = \inf \{d(x,t)\colon (x,t) \in A\}$, and suppose that 
$\gamma < d(x,t)$, $(x,t) \in A $. Then there exists a sequence $((x_{n},t_{n}))_{n \in \mathbb{N}}$ 
such that $(d(x_{n},t_{n}))_{n \in \mathbb{N}}$ decreases to $\gamma$. For each $m < n$ we have
\begin{equation*}
 \delta(\xm,\xn) \leq \delta(\xm,x_{m+1})+ \cdots + \delta(x_{n-1},\xn) \leq d(\xm,t_{m}) - d(\xn,\tn),
\end{equation*}
and therefore, $\xnn$ is a Cauchy sequence. Similarly, $(t_{n})_{n \in \mathbb{N}}$ is also a Cauchy sequence. 
Then for $x = \li \xn$, $t = \li \tn$, $(x,t) \in G$ holds as $G$ is closed. For any 
$(z,v) \in A$, and large $n$ we have
\[
 \delxz \leq \delta(x,x_{n}) + \delta(x_{n},z) \leq  \delta(x,x_{n}) + d(z,v) - d(\xn,\tn).
\]
Both metrics are equivalent on $X \cup F(X)$, and consequently, $\delxz \leq d(z,v) - d(x,t)$. 
In a similar way we state that $\delta(t,v) \leq k (d(z,v) - d(x,t))$. Therefore, $(x,t) \in A$, and 
$d(x,t) = \gamma$. Suppose $x \notin F(x)$, and let $(z,v) \neq (x,t)$ be as in \eqref{co16}, {\it i.e.} 
\[
 d(z,v) \leq d(x,t) - \delxz < d(x,t) = \gamma.
\]
On the other hand, for any $(y,s) \in A$ we have  
\begin{equation*}
\begin{split}
 &\delta(y,z) \leq \delta(y,x) + \delxz \leq d(y,s) - d(x,t) + d(x,t) - d(z,v) = \\
 &d(y,s) - d(z,v) \mbox{,  } \qquad \delta(s,v) \leq k (d(y,s) - d(z,v)).
\end{split}
\end{equation*}
Consequently, $(z,v) \in A$ holds, and therefore, $\gamma \leq d(z,v)$ which contradicts 
$d(z,v) < \gamma$.  $\quad \square$
\par Now, we are ready to prove an analog of Theorem \ref{Th34}.
\begin{theorem}
\label{Th39}
Let $X$ be a closed set in a complete metric space $(Y,d)$. If $ \FXY $ is a closed valued 
$\alpha$-contraction, $\aeps < 1$, and
\begin{equation}
\label{co17}
\begin{split}
 &\mbox{for an  } \ej \in (0,\epsilon) \mbox{, and each  } \xXFx \mbox{, } t \in F(x) 
     \mbox{ there exists} \\
 &\mbox{a  } \zXx \mbox{ such that } (1-\ej) d(x,z) \leq d(x,t) - d(z,t) 
\end{split}
\end{equation}
holds, then the following is satisfied:
\begin{equation}
\label{co18}
\begin{split}
 &\mbox{for each } \xXFx \mbox{, } t \in F(x) \mbox{ there exist } \zXx \mbox{, } v \in F(x) \\
 &\mbox{such that } (1-\alpha-\epsilon) \max\{d(x,z),d(t,v)/(\aeps)\} \leq d(x,t) - d(z,v), 
\end{split}
\end{equation}
and $F$ has a fixed point.
\end{theorem}
Proof. 
For each $\ej \in (0,\epsilon)$, $x,z \in X$, and $t \in F(x)$ there exists a $v \in F(z)$
such that
\begin{equation}
\label{co19}
 d(t,v) \leq (\aeps - \ej) d(x,z) < (\aeps) d(x,z)
\end{equation}
holds, as $F$ is an $\alpha$-contraction. On the other hand, in view of \eqref{co17}, and \eqref{co19}
we have
\begin{equation*}
\begin{split}
 &d(z,v) + (1 - \ej)d(x,z) \leq d(z,v) + d(x,t) - d(z,t) \leq \\
 &d(t,v) + d(x,t) \leq (\aeps - \ej) d(x,z) + d(x,t).
\end{split}
\end{equation*}
Hence it follows that 
\[
 (1-\alpha-\epsilon) d(x,z) \leq d(x,t) - d(z,v).
\]
Now, it is clear that \eqref{co16} is satisfied for $\delta = (1-\alpha-\epsilon)d $, $k = \aeps$ ({\it i.e.} 
\eqref{co18} is valid), and we may apply Theorem \ref{Th38}, as the graph of a closed valued $\alpha$-contraction 
on a closed set is closed. $\quad \square$
\par From Lemma \ref{Le35}, and Theorem \ref{Th39} we obtain the following refinement of Maciejewski's 
Theorem 2.1 \cite[p. 1210]{Ma}, and of the well known theorem of Lim for inward contractions 
\cite[Theorem 1]{Lw}:
\begin{theorem}
\label{Th310}
Let $X$ be a closed set in a complete metric space $(Y,d)$. If $ \FXY $ is a closed valued 
$\alpha$-contraction such that $F(x) \subset \tIXx$, $x \in X$, then for each $\epsilon > 0$
with $\aeps < 1$ condition \eqref{co18} is satisfied, and $F$ has a fixed point.
\end{theorem}
\par The proofs of Theorems \ref{Th38}, \ref{Th39} yield
\begin{corollary}
\label{Cor311}
Let $X$ be a closed set in a complete metric space $(Y,d)$. If $ \FXY $ is a closed valued 
$\alpha$-contraction such that $F(x) \subset \tIXx$, $x \in X$, then for each $\epsilon > 0$
with $\aeps < 1$ there exist convergent sequences $\xnn $, $(t_{n})_{n \in \mathbb{N}}$ such that
$\tn \in F(\xn)$, $n \in \mN$, $x = \li \xn \in F(x)$, $t = \li \tn \in F(x)$, and
\[
(1-\alpha-\epsilon) d(\xn,\xnp) \leq d(\xn,\tn) - d(\xnp,t_{n+1}), \quad n \in \mN.
\]
\end{corollary}
\par Theorem \ref{Th39} is also a consequence of the next ``tool'' theorem.
\begin{theorem}
\label{Th312}
Let $X$ be a set in a complete metric space $(Y,d)$. Assume that $\FXY$ is a mapping, and
$\delta$ is a complete metric on $G = graph(F)$. If  
\begin{equation}
\label{co20}
\begin{split}
 &\mbox{for each  } \xXFx \mbox{, } t \in F(x) \mbox{ there exist } \zXx \mbox{, }  \\
 &v \in F(x) \mbox{ such that } \delta((x,t),(z,v)) \leq d(x,t) - d(z,v)
\end{split}
\end{equation}
holds, then $F$ has a fixed point.
\end{theorem}
Proof. 
Once again we use the method from the proof of Theorem \ref{Th32}. For simplicity of notations let us
adopt $p=(x,t)$, $\pn=(\xn,\tn)$, $q=(z,v)$, and (not quite correctly) $d(p) = d(x,t)$,
$d(\pn) = d(\xn,\tn)$, $d(q) = d(z,v)$. If $\xo \notin F(\xo)$, then there exists a maximal 
set $A \subset G$, $A \neq \emptyset$ such that all $p,q \in A$ satisfy
\[
 \delta(p,q) \leq \mid d(p) - d(q) \mid.
\]
Let us adopt $\gamma = \inf \{ d(q)\colon q \in A\}$, and suppose that $\gamma < d(q)$, for each 
$q \in A$. Then there exists a sequence $(p_{n})_{n \in \mathbb{N}}$ in $G$ such that 
$(\varphi(p_{n}))_{n \in \mathbb{N}}$ decreases to $\gamma$. Hence it follows that 
\[
 \delta(p_{m},\pn) \leq d(p_{m}) - d(\pn), \quad m < n,
\]
{\it i.e.} $(p_{n})_{n \in \mathbb{N}}$ is a Cauchy sequence, and there exists a $p = \li \pn \in G$ as $\delta$ is 
complete. For any $q \in A$, and large $n$ we have
\[
 \delta(p,q) \leq \delta(p,\pn) + \delta(\pn,q) \leq \delta(p,\pn) + d(q) - d(\pn),
\]
and $\delta(p,q) \leq d(q) - d(p)$, which means $p \in A$, and $d(p) = \gamma$ (implies 
$p$ is unique in $A$). Suppose $x \notin F(x)$. Then for $q =(z,v)$ as in \eqref{co20} from 
$z \neq x$ it follows that $d(q) < \gamma$, and for any $r=(y,s) \in A$ we obtain
\[
 \delta(r,q) \leq \delta(r,p) + \delta(p,q) \leq d(r) - d(p) + d(p) - d(q) = d(r) - d(q),
\]
which means $q \in A$, and $\gamma \leq d(q)$ - a contradiction.  $\quad \square$
\par The previous theorem could be used to prove Theorem \ref{Th310}. The respective complete metric $\delta$
on $graph(F)$ is defined by the left side of the inequality in condition \eqref{co18}.
\par Let us go back to condition \eqref{co15}, and Theorem \ref{Th36}. As regards the ``nonexpansive'' case 
($ \aeps = 1 $), \eqref{co15} means
\begin{equation*}
\begin{split}
 &\mbox{for each } x \in X \setminus F(x) \mbox{ there exists a } z \in X \setminus \{x\}  \\
 &\mbox{such that }  \dzF  \leq \dxF, 
\end{split}
\end{equation*}
and it would be rather an ambitious task to prove anything for such a condition ;). More interesting is
its sharper form:
\begin{equation}
\label{co21}
\begin{split}
 &\mbox{for each } x \in X \setminus F(x) \mbox{ there exists a } z \in X \\
 &\mbox{such that } \dzF < \dxF.
\end{split}
\end{equation}
\par It is known (see \cite[proof of Theorem 8]{Lo}) that if $X$ is a bounded closed convex set in
a Banach space $Y$, and $\FXY$ is a nonexpansive weakly inward mapping, then there exists a $\xnn$ 
such that $\li d(\xn,F(\xn)) = 0$. Consequently, for such a case, and $\dxF = 0$ iff 
$x \in F(x)$ \eqref{co21} holds. 
Therefore, the following is also an ``inward'' theorem ({\it cp.} Theorems \ref{Th15}, \ref{Th17}).
\begin{theorem}
\label{Th313}
Let $X$ be a compact set in a metric space $\Yd$, and let $\FXY$ be a mapping satisfying \eqref{co21}, 
and such that $\dFc$ is lower semicontinuous. Then $F$ has a fixed point.
\end{theorem}
Proof. 
From the lower semicontinuity of $\dcF$ it follows that for each $\lambda \geq 0$ the set
$X_{\lambda} = \{x \in X \colon \dxF \leq \lambda \}$ is compact. Suppose 
$\alpha = inf \{\lambda \colon X_{\lambda} \neq \emptyset \} > 0$. $X_{\alpha}$ is clearly nonempty, as it is an
intersection of a decreasing family of nonempty compact sets. For any $x \in X_{\alpha}$ there
exists a $z \in X$ such that $\dzF < \dxF$ (see \eqref{co21}) - a contradiction.
Thus $X_{0}$ consisting of fixed points of $F$ is nonempty. $\quad \square$ 

\section*{Acknowledgements}
This work was partially supported by the Faculty of Applied Mathematics AGH UST statutory tasks within subsidy 
of the Polish Ministry of Science and Higher Education, grant no. 16.16.420.054.
\par 
\pagebreak
\mbox{Faculty of Applied Mathematics} \linebreak
\mbox{AGH University of Science and Technology} \linebreak
\mbox{Al. Mickiewicza 30} \linebreak
\mbox{30-059 KRAK\'OW, POLAND} \linebreak
\mbox{E-mail: pasicki@agh.edu.pl}


\begin{thebibliography}{00}

\bibitem{Br}
F.E. Browder, The fixed point theory of multi-valued mappings in topological
vector spaces, Math. Ann., 177 (1968) 283-301. 
\bibitem{Ca}
J. Caristi, Fixed point theorems for mappings satisfying inwardness conditions, Trans. Amer. Math. Soc.,
215 (1976) 241-251.
\bibitem{Do}
D. Downing, W.A. Kirk, Fixed point theorems for set-valued mappings in metric and Banach spaces,
Math. Japonica, 22 (1977) 99-112. 
\bibitem{Ed}
M. Edelstein, The construction of an asymptotic center with a fixed-point property,
Bull. Amer. Math. Soc., 78 (1972) 206-208. 
\bibitem{Go}
K. Goebel, On a fixed point theorem for multivalued nonexpansive mappings,
Ann. Univ. Mariae Curie-Sk{\l}odowska Sect. A, 29 (1975) 69-71. 
\bibitem{Lf}
T.C. Lim, A fixed point theorem for multivalued nonexpansive mappings in uniformly
convex Banach space, Bull. Amer. Math. Soc., 80 (1974) 1123-1126.
\bibitem{Lo}
T.C. Lim, On asymptotic centers and fixed points of nonexpansive mappings, Canad. J. Math.,
32 (1980) 421-430.
\bibitem{Lw}
T.C. Lim, A fixed point theorem for weakly inward multivalued contractions, J. Math. Anal. Appl.,
247 (2000) 323-327. 
\bibitem{Ma}
M. Maciejewski, Inward contractions on metric spaces, J. Math. Anal. Appl., 330 (2007) 1207-1219. 
\bibitem{Na}
S.B. Nadler, Multi-valued contraction mappings, Pacific J. Math., 30 (1969) 475-488.
\bibitem{Pc}
L. Pasicki, A short proof of the Caristi theorem, Comment. Math., 20 (1978) 427-428. 
\bibitem{Pf}
L. Pasicki, A basic fixed point theorem, Bull. Pol. Acad. Sci. Math., 54 (2006) 85-88.
\bibitem{Ps}
L. Pasicki, Bead spaces and fixed point theorems, Topology Appl., 156 (2009) 1811-1816,
DOI: 10.1016/j.topol.2009.03.042. 
\bibitem{Pt}
L. Pasicki, Towards Lim, Topology Appl., 158 (2011) 479-483, DOI: 10.1016/j.topol.2010.11.024.  
\bibitem{Pu}
L. Pasicki, Uniformly convex spaces, bead spaces, and equivalence conditions,
Czechoslovak Math. J., 61 (2011) 383-388.
\bibitem{Pv}
L. Pasicki, Transitivity and variational principles, Nonlinear Anal., 74 (2011), 5678-5684,
DOI: 10.1016/j.na.2011.05.054. 
\bibitem{p30}
L. Pasicki, Variational principles and fixed point theorems, Topology Appl., 159 (2012), 3243-3249, 
DOI: 10.1016/j.topol.2012.07.002.
\bibitem{Xu}
H.K. Xu, Multivalued nonexpansive mappings in Banach spaces, Nonlinear Anal.,
43 (2001) 693-706. 
\end{thebibliography}
\end{document}